\documentclass[a4paper,12pt]{article}
\usepackage{amsmath}
\usepackage{amssymb, amsbsy, amstext}
\usepackage{srcltx}
\usepackage{amscd,amsfonts,color}
\input amssym.def
\input amssym.tex

\newtheorem{lem}{Lemma}[section]%
\newtheorem{theorem}[lem]{Theorem}%
\newtheorem{defi}[lem]{Definition}%
\newtheorem{exam}[lem]{Example}%
\newtheorem{prop}[lem]{Proposition}%

\def\a{\alpha} \def\b{\beta}   
 \def\s{\sigma} \def\t{\tau}

 \def\O{\Omega}

    \def\oc{\overline C}
\def\ob{\overline B}

 \def\ox{\overline X}  \def\o1{\overline 1}

\def\ola{\overline a}    
\def\olb{\overline b}

\def\olz{\overline z}  

\def\o{\overline}   \def\olb{\overline b}
\def\di{\bigm|} \def\lg{\langle} \def\rg{\rangle}

\def\Aut{\hbox{\rm Aut\,}} \def\Inn{\hbox{\rm Inn}} \def\Syl{\hbox{\rm Syl}}
  
\def\Cay{\hbox{\rm Cay }}  \def\mod{\hbox{\rm mod }}
  
 \def\PG{\hbox{\rm PG}} \def\PGL{\hbox{\rm PGL}}
 \def\AGL{\hbox{\rm AGL}} \def\GL{\hbox{\rm GL}} \def\bGL{\hbox{\bf{\rm GL}}}  \def\P\GL{\hbox{\rm P\GL}}
  \def\FF{{\hbox{\sf F\kern-.43emF}}}
\def\cal{\mathcal} \def\char{\hbox{\rm char}}

\def\Z{\hbox{\rm Z}}

\def\o{\hbox{\rm o}}
\def\ord{\hbox{\rm ord}}
\def\lcm{\hbox{\rm lcm}}
\def\char{ \, {\rm char}\,}

\def\calm{\mathcal{M}} \def\MM{\mathcal{calm}}

  \def\MF{{\hbox{\sf MF\kern-.43emF}}}
\def\ZZ{\mathbb{Z}}  
\def\CM{\mathbb{CM}}

\def\nd{\mathrel{\bigm|\kern-.7em/}} 
 \def\f{\noindent}
\def\qed{\hfill $\Box$} \def\demo{\f {\bf Proof}\hskip10pt}

  \def\GG{{\cal G}}  \def\MM{{\cal M}} 

\def\bF{\hbox{\bf F}}

\begin{document}

\begin{center}
{\bf\large  Regular Cayley Maps of Elementary abelian $p$-groups:  Classification and Enumeration}
\end{center}

%\bigskip

\begin{center}
Shaofei Du,  Hao Yu{\small \footnotemark} and Wenjuan Luo\\
\medskip
{\it {\small
Capital Normal University,\\ School of Mathematical Sciences,\\
Beijing 100048, People's Republic of China
}}
\end{center}

\footnotetext{Corresponding author: 3485676673@qq.com. This work is supported in part  by the National Natural Science Foundation of China (12071312 and 11971248).}

\renewcommand{\thefootnote}{\empty}%{footnote}}

\begin{abstract}
Recently, regular Cayley maps of cyclic groups and dihedral groups have been classified in \cite{CT} and \cite{KK1}, respectively. A nature question is to classify  regular Cayley maps of
elementary abelian $p$-groups $\ZZ_p^n$.  In this paper, a complete  classification of regular Cayley maps of $\ZZ_p^n$ is given and moreover, the number of these maps and their genera  are enumerated.
  \end{abstract}

\section{Introduction}
A (topological) map is a cellular decomposition of a closed surface. A common way to describe maps is to view them as 2-cell embeddings of graphs. An orientation preserving automorphism of an orientable  map is an automorphism of the underlying graph which extends to an orientation preserving self-homeomorphism of the supporting surface.  The group of all such automorphisms  acts semiregularly on the  arcs of the underlying graph and in an extreme case, when the action is regular, the map itself is called regular. A regular embedding of a graph is a 2-cell embedding of a graph into a surface in a way that the associated map is regular. Regular maps have been studied in connection with various branches of mathematics including Riemann surfaces and algebraic curves. For more information about regular maps and their connections to other fields of mathematics we refer the reader to \cite{Jon,Jon1,JS1,JS2,NS1,RSJ,SS2}.

\vskip 3mm

Given a group $G$ and a generating subset $S$ of $G$ such that $S=S^{-1}$ and $1\not\in S$, a {\it Cayley graph} $\Cay(G,S)$ is the graph with vertices $G$ and arc-set $\{ (g, gs)\di g\in G, s\in S\}$.
Take a permutation $\rho $ on $S$,  a Cayley map $\CM:=\mathrm{Cay}(G,S,\rho)$ is an embedding
of $\Cay(G,S)$ into an orientable closed surface such that, at each vertex $g\in G$, the local
orientation $R_g$ of the darts $(g,gs)$ incident with $g$ agrees with a prescribed
cyclic permutation $\rho$ of the generating set $S$, that is, $R_g(g,gs)=(g,g\rho(s))$ for all
$g\in G$ and $s\in S$. Clearly, the (orientation preserving) automorphism group $\Aut(\CM)$ of a Cayley map $\CM$ contains
the  left regular representation  of $G$, which acts transitively on vertices.  If the
(cyclic) stabiliser of a vertex is transitive on its adjacent vertices, then the automorphisms
group $\Aut(\CM)$ is regular on the arcs. In this case $\CM$ is called a
\textit{regular Cayley map}.

\vskip 3mm
The term {\it Cayley maps} appeared firstly in a paper of Biggs \cite{Big} in 1972, which are now known as {\it balanced Cayley maps}, where
a Cayley map $\CM(G, S, \rho)$  is called {\it balanced}  if $\rho(x^{-1})=\rho(x)^{-1}$, for any $x\in S$. It was shown  in \cite[Proposition 2.3]{CJT}
that a regular Cayley map for the group $G$ is balanced if and only if the
associated skew morphism is an automorphism.
An  important example is  the solution of the famous Heawood map coloring problem
(see \cite{RM}):
the cases $n \equiv 0, 4 ~\rm{and}~ 7 {\pmod {12}}$ of the solution rely
on constructing a triangular Cayley map of the complete graph of order $n$ (see also \cite{GT}).
Another example  is the classification of regular maps of complete graphs; by a result in \cite{JJ},
all such maps are balanced Cayley maps. Cayley maps were systematically  studied in \cite{SS2,SS1}, mainly  on balanced and antibalanced Cayley maps.
In general, map automorphism groups of arbitrary Cayley maps have been described in \cite{Jaj1,Jaj2}.
For other results in Cayley  maps, the reader may refer  \cite{AR,GT,Tuc}.

\vskip 3mm
As an algebraic tool to investigate regular Cayley maps, {\it skew morphisms} of a group were introduced by Jajcay and \v{S}ir\'a\v{n} in  \cite{JS}, while for the definition, see next section. They showed that a Cayley map $\CM$ is regular if and only if $\rho$ extends to a skew-morphism of $G$. Moreover, if $\CM$ is regular, then $\Aut(\CM) = G\lg\s\rg$, where $\s$ is a skew-morphism of $G$ (see \cite[Theorem 1.1]{JS}). For  partial results  on   skew-morphisms of cyclic groups, see  ~\cite{CJT,CT,DH,Kwo,KN1,KN2}; that of  dihedral groups see \cite{CJT,KK2,KMM,Zhang2,ZD}; that of  monolithic groups (in particular, finite nonabelian simple groups),  see  see \cite{BCV}; that of   finite nonabelian characteristically    simple groups,  see ~\cite{CDL}; and that of
 elementary abelian $p$-groups $\ZZ_p^n$, see \cite{DYL}.

\vskip 3mm Conder and  Tucker  gave a complete classification of all regular Cayley maps for finite cyclic groups in \cite{CT}. Recently,  Kov\'acs and Kwon  finished  a complete classification of regular Cayley maps of dihedral groups in \cite{KK1}.  Now a natural problem  is posed:
 \vskip 3mm{\it  Classify regular Cayley maps for elementary abelian $p$-groups $\ZZ_p^n$. Is there any unbalanced  regular Cayley maps of $\Z_p^n$ ?}
 \vskip 3mm
 In this paper, this problem will be   completely solved,  containing the following five   theorems. The first one  gives a characterization of the automorphism group of  a   regular Cayley map of $G\cong \ZZ_p^n$.
 \begin{theorem}  \label{main1} Let $\CM$ be a regular Cayley map of $G\cong \ZZ_p^n$ such that   a point-stabiliser $\lg \s\rg $ of $\Aut(\CM)$ is of order $kp^m$ where $m\ge 0$ and $p\nmid k$.
 Then  $\Aut(\CM)=(G\rtimes  \lg \s^k\rg )\rtimes \lg \s^{p^m}\rg$. Moreover, $\Aut(\CM)$ contains a normal vertex regular subgroup $T\cong \ZZ_p^n$ such that either $T=G$ or $T\cap G\cong \Z_p^{n-1}$, where $n\ge 2$.
  \end{theorem}

In this paper, as usual, an orientably-regular map  with the automorphism group $X$  will be represented by  $\MM(X; r,\ell)$,  where  $r$ corresponds to the local rotation and $\ell$ the arc-revising, so called  {\it algebraic map}.
Moreover, from  Theorem~\ref{main1}, $\Aut(\CM)=T\rtimes \lg \s \rg\le \AGL(n,p)$.
So we identify $T$ with an $n$-dimensional row space $V(n,p)$ and  $\s$ with an element in $\bGL(n,p)$
so that $\s^{-1}t\s=\s(t)$, for any $t\in T$.
To state next theorem, we introduce two subsets of $\bGL(n,p)$:
\vskip 4mm {\it Let  ${\bf M}(n,p)$ be the set of representatives of conjugacy classes of  linearly transformation  $\s$ in $\bGL(n,p)$ whose  minimal polynomial  is of degree $n$ and $-I \in \lg \s\rg$, where $I$ denotes  the identity, noting that
   $-I=I$ if $p=2$. Let  ${\bf M}_1(n,p)$ is  the set of  elements $\s$ in ${\bf M}(n,p)$    such that  $p\di o(\s)$ and $\s^{\frac{o(\s)}p}$
    fixes an $(n-1)$-dimensional subspace  pointwise.}

  \vskip 3mm Then we state the second  main theorem, noting that in the theorem, $t(-I)$ is the product of $t\in T$ and $-I\in \lg \s\rg $ in the affine subgroup $T\rtimes \lg \s\rg $.

\begin{theorem} \label{main2} Every regular  Cayley map  of $G\cong \ZZ_p^n$  is isomorphic to
$\MM(T\rtimes \lg \s\rg ;\, \s,\, t(-I))$,
where $\s \in {\bf M}(n,p)$ and $t$ may be taken  from  $T$ such that     $\{\s^i(t)\di 0\le i\le n-1\}$ is a base for  $T$.
Moreover, given $p$ and $n$, different  $\s $ in   ${\bf M}(n,p)$ give nonisomorphic   maps.
\end{theorem}

Remark that
\vskip 3mm

(1)  Under the base $\{\s^i(t)\di 0\le i\le n-1\}$, ${\bf M}(n,p)$  may be identified with the  set of   matrices of the form  $$\left(\small\begin{array}{ccccc}
    0  & 0 & \cdots & 0 & {r_0} \\
    1  & 0 & \cdots & 0 & {r_1} \\
    & & \cdots& & \\
     0 & 0&\cdots &1& {r_{n-1}}\\
    \end{array}\right)\in{\bf M}(n,p),$$
where $r_0, r_1,\cdots, r_{n-1}\in \bF_p$ and $r_0\ne 0$    and $t$ may be chosen as $(1, 0, \cdots, 0)$.  In particular, two such different matrices   are not conjugate each other.

\vskip 3mm
(2)  ${\bf M}_1(n,p)$ may be an empty set for some $n$ and $p$.

\vskip 3mm
(3)    In Theorem~\ref{main1},   if    $T\cap G\cong \Z_p^{n-1}$  for some $\s$, then the map
$\MM(T\rtimes \lg \s \rg ;\, \s, \, t(-I))$  is unbalanced for $G$.

\begin{theorem} \label{main3}
The set ${\bf M}(n,p)$ and ${\bf M_1}(n,p)$ are characterized in Lemmas~\ref{Mnp}.
Moreover,  the case  $T\cap G\cong \Z_p^{n-1}$ in Theorem~\ref{main1} happens   if and only if $p\ne 2$ and  $\s \in {\bf M}_1(n,p)$ when ${\bf M}_1(n,p)\ne \emptyset $.  \end{theorem}

\begin{theorem} \label{main4}
The number of elements  in  ${\bf M}(n,p)$ is enumerated in Lemma~\ref{nMnp}.
\end{theorem}

\begin{theorem} \label{main5}
The genera of regular Cayley maps $\MM(T\rtimes \lg \s\rg ;\, \s,\, t(-I))$  are computed in Lemma~\ref{genus}, where $\s\in {\bf M}(n,p).$
\end{theorem}

After this introductory section, some  preliminary results will be given in Section 2, Theorems~\ref{main1}-\ref{main5} will be proved in Sections 3-7, respectively.
To end up this section, we list all notations which will be used in this paper:
\vskip 3mm
{\it
$\gcd(a_1,a_2)$: the greatest common factor of integers $a_1$ and $a_2$;

$\lcm(a_1, a_2, \cdots , a_n)$ and $\lcm\{a_i\mid 1\leq i\leq n\}$: the least  common multiple of integers

 $a_1, a_2,\cdots,a_n$;

$(f(x), g(x))$: the greatest common factor of polynomials $f(x)$ and $g(x)$;

${\bf F}_q$ and ${\bf F}_q^*$: the  field of order $q$ and its multiplicative group, resp.;

$V(n,p)$: $n$-dimensional row space over ${\bf F}_q$;

$PG(V)$: the projective space of $V$;

$[u]_p$: the maximal integer $k$ such that $p^k\di u$;

${\bf F}_q[x]$  and ${\bf MF}_q[x]$:  polynomial ring over $\bF_q$ and    monic polynomials in it, resp. ;

 $I$: identity transformation;

 $(a_{ij})_{n\times n}$: $n$-th degree matrix;

 $|M|$: the cardinality of a set $M$;

 $|G|$,  $o(g)$:   the order of a group $G$ and an element $g\in G$, resp.;

 $H\lhd G$:  $H$ is a proper normal subgroup of $G$;

$H\char~G$: $H$ is a characteristic subgroup of $G$;

 $G'$ and $Z(G)$: the derived subgroup and the center of $G$, resp.;

 $G_X:=\cap_{x\in X}G^g$:  the core of $G$ in $X$;

 $\Syl_p(G)$:  the set of Sylow $p$-subgroups of $G$;

 $\O_1(P)$: the group generated by all elements of order $p$ in a $p$-group $P$;

 $[A, B]:$ the subgroup generated by commutators $[a,b]:=a^{-1}b^{-1}ab$, $a\in A$ and $b\in B$;

 $M\rtimes N:$   a semidirect product of $M$ by $N$, in which $M$ is  normal.}
\vskip 3mm

\section{Preliminary Results}\label{2}
In the investigation of maps, it is often useful to replace topological maps on surfaces with their combinatorial counterparts. Let $\GG=\GG(V,D)$ be a  graph with vertex set $V=V(\GG)$ and arc (dart) set $D=D(\GG)$.  It is well-known that graph embeddings into orientable surfaces can be described by means of local rotations (see \cite{GT,JS2}). A map $\MM$ with underlying graph $\GG$ can be identified with a triple ${\MM}=\MM (\GG;R,L)$, where $R$ is a rotation and $L$
 is the arc-reversing involution of $\GG$.
Given a group $X=\lg r,\ell\rg$ where $\ell^2=1$, from the above arguments, we may define an {\it algebraic map} $\MM(X;r,\ell)$ as follows: set $D=G$ and consider the left multiplication action $L(G)$ of $G$ on $D$. The vertices, edges and faces are just the orbits of $\lg r\rg$, $\lg \ell\rg$ and $\lg r\ell\rg$, respectively, with the natural inclusions relation.
Moreover,  two  maps $\MM(X_1; r_1, \ell _1)\cong \MM(X_2; r_2, \ell _2)$ if and only if there
 exists an group isomorphism $f $  from $X_1$ to $X_2$  such that $f(r_1)=r_2$ and $f(\ell _1)=\ell _2$.

\vskip 3mm
A \textit{skew-morphism} of a finite  group $G$ is a permutation $\s$ on $G$ fixing the
identity element, and for which there exists an integer function $\pi$ on $G$ such that
$\s(gh)=\s(g)\s^{\pi(g)}(h)$ for all $g,h\in G$.  Given a skew-morphism $\pi$ of $G$,  $X:=L_{G}\langle\s \rangle$
is a permutation group on $G$, called the \textit{skew-product group} of $\s$,  see
\cite{CJT, ZD}.  Sometimes, simply we just say skew-product groups of $G$.

\vskip 3mm
It was shown by Jajcay and \v{S}ir\'{a}\v{n} that a Cayley map $\calm$ is
regular if and only if $\rho$ extends to a skew-morphism of $G$, see \cite[Theorem~1]{JS}.
Thus the problem
of determining all regular Cayley maps of a group $G$ is equivalent to the problem of determining
all skew-morphisms of $G$ containing a generating orbit which is closed under taking inverses.
In \cite{DYL},  the  authors characterized   the skew-product groups of all skew-morphism of elementary abelian $p$-groups.
This is the start point of the present paper.

\begin{prop}\label{dty1}\cite[Theorem 1.1]{DYL}
  Let $G\cong \ZZ_p^n$ where $p$ is a prime.   Let $X=G\lg \s \rg $ be  the  skew-product group of  a skew-morphism $\s$ of $G$, with order $|\s|=k\cdot p^m\ge 2$, where $m\ge 0$ and $p\nmid k$.
   Let $P=G\lg \s^k\rg $.
  Then  $P\lhd X$, either $G\lhd X$ or $G=G_X\times B$ where $B\cong \ZZ_p$. More precisely,   one of the following holds:
\begin{enumerate}
  \item[\rm(1)] Either $m=0$ or $p=2$:  $X=G\rtimes \lg \s\rg $;
  \item[\rm(2)] $p\ne 2$, $m\ge 1$, $k=1$:  either $X=G\rtimes \lg \s\rg $; or  $X=(G_X\rtimes \lg \s \rg)\rtimes B$;
  \item[\rm(3)]  $p\ne 2$, $m\ge 1$, $k\ge 2$: either  $X=G\rtimes \lg \s\rg $; or

   $G\lhd P$,  $X=((G_X\times B)\rtimes \lg \s^k\rg )\rtimes \lg \s^{p^m}\rg $; or

   $G\ntriangleleft P$, $X=((G_X\rtimes \lg \s^k\rg)\rtimes B)\rtimes \lg\s^{p^m}\rg $.
\end{enumerate}
\end{prop}
\vskip 3mm
Several results on  finite fields $\bF_q$ of order $q$ with the character $p$  will be needed.

\begin{prop} \label{RN4}\cite[Theorem 2.14]{RN}
If $f(x)$ is an irreducible polynomial in $\bF_q[x]$ of degree $m$, then $f(x)$ has a root $\a$ in $\bF_{q^m}$. Furthermore, all the roots of f are simple and are given by the number  $m$ distinct elements
$\a, \a^q,\cdots, \a^{q^{m-1}}$ of $\bF_{q^m}.$\end{prop}

\begin{defi} \label{RN1}\cite[Definition 3.2]{RN}
Let $f(x)\in \bF_q[x]$ be a nonzero polynomial. If $f(0)\neq0$, then the least positive integer $e$ for which $f(x)$ divides $x^e - 1$ is called the order of $f(x)$, denoted by $\ord(f(x))$ .
\end{defi}

\begin{prop} \label{RN2}\cite[Theorem 3.9]{RN}
Let $f(x)=\Pi_{i=1}^rp_i^{k_i}(x)$, where $p_1(x), p_2(x) \cdots, p_r(x)$ are monic  pairwise relatively prime   irreducible  polynomials over $\bF_q$.
Let $t_i$  be the smallest integer with $p^{t_i}\geq k_i$. Then
 $$\ord(f(x))=\lcm(\ord(p_1(x)p^{t_1},\cdots ,  \ord(p_2(x)p^{t_2}, \ord(p_r(x)p^{t_r}).$$
\end{prop}

\begin{prop}\label{RN3}\cite[Theorem 3.5]{RN}
The number of monic irreducible polynomials in $\bF_q[x]$ of degree $m$ and order $e$ is equal to $\frac{\varphi(e)}m$ if $e\geq2$ and $m$ is the multiplicative order of $q$ modulo $e$, equal to $2$ if $m = e =1$, and equal to $0$ in all other cases. In particular, the degree of an irreducible polynomial in $\bF_q[x]$ of order $e$ must be equal to the multiplicative order of $q$ modulo $e$.
\end{prop}

The last lemma deals with isomorphisms between two affine groups.

\begin{lem} \label{dt1}
Let $X_1=T\rtimes \lg \s_1\rg$ and $X_2=T\rtimes \lg\s_2\rg$ be two subgroups of $\AGL(n,p)$, where $T$ is the translation subgroup, and $\s_1$ and $\s_2$ are nontrivial linearly transformations  in $\GL(n,p)$.
Suppose $f$ is an isomorphism from $X_1$ to $X_2$ fixing $T$ setwise and mapping $\s_1$ to $\s_2$. Then there exists an element $u\in \GL(n,p)$ such that $f=Inn(u)|_{X_1\rightarrow X_2}$, where $Inn(u)$ is the inner automorphism
of $\AGL(n,p)$ induced by $u$. In particular, if $\lg\s_1\rg=\lg\s_2\rg$, which is fixed by $f$,  then $u\in N_{\GL(n,p)}(\lg\s_1\rg)$.
\end{lem}
\demo The lemma can be equivalently proved  by a matrix representation. Take $T=\lg t_1,t_2,\cdots,t_n\rg$ where
$t_i=t_{(0,\cdots,1,\cdots,0)}$, where 1 is in the $i$th-component.
Let $f(t_i) = t_1^{u_{i1}} t_2^{u_{i2}}\cdots t_n^{u_{in}}$ for some $u_{ij}\in\ZZ_p$.
Then $f(t_i)= t_{(u_{i1},u_{i2},\cdots,u_{in})}$.
Moreover, it is easy to check that for any $\a\in V(n,p)$, $f(t_{\a})=t_{\a u}=t_{\a}^u$
for $u=(u_{ij})_{n\times n}\in \GL(n,p)$.
Because $f(\s^{-1}t_{\a}\s) = f(t_{\a}^\s)$ and $f(\s)\in \GL(n,p)$,
we have $t_{\a u f(\s)}=t_{\a\s u}$, which implies that
$f(\s)=u^{-1}\s u$. Therefore $f=Inn(u)|_{X_1\rightarrow X_2}$. The second part of the lemma is immediate.
\qed

\section{Structures of $\Aut(\CM)$}
Let $\CM$ be a regular Cayley map of $G\cong \ZZ_p^n$ and $X=\Aut(\CM)$, with
 a point-stabiliser $C=\lg\s\rg\cong\ZZ_{kp^m}$, where $p\nmid k$.
Then
 $X=GC$, where the core $C_X$ is trivial  and  $o(\s)$ is even if $p\ne 2$.
   Moreover, set $C_1=\lg\s_1\rg$ and $C_2=\lg\s_2\rg$, where $\s_1:=\s^k$ and $\s_2:=\s^{p^m}$.
   If $m\geq 1$, set $z:=\s^{kp^{m-1}}$ which is contained in $C_1.$
   Then  $P=GC \in \Syl_p(X)$.   Moreover, $X=\lg \s, \ell\rg $, where $\ell$ is an involution.
 In this section, we shall prove   Theorem~\ref{main1}.
 \vskip 3mm
\f {\bf Theorem 1.1} {\it  $\Aut(\CM)=(G\rtimes  \lg \s^k\rg )\rtimes \lg \s^{p^m}\rg$ and   $\Aut(\CM)$ contains a normal vertex regular subgroup $T\cong \ZZ_p^n$ such that either $T=G$ or $T\cap G\cong \Z_p^n$ where $n\ge 2$.}
\vskip 3mm

 If $p=2$, then $X=G\rtimes C$ by Proposition~\ref{dty1} and so the theorem holds by taking $T=G$. Hence in the remaining of this section,  we assume $p$ is an odd prime and
 carry out the proof  by the following Lemmas~\ref{normal1} and \ref{normal4}.

\begin{lem} \label{normal1}
$G\lhd P$ and $X=(G\rtimes C_1)\rtimes C_2$.
\end{lem}
\demo  By Proposition~\ref{dty1}, $P\lhd X$. Suppose that $G\lhd P$. Then it follows from  Proposition~\ref{dty1} that  $X=(G\rtimes  C_1)\rtimes C_2,$  as desired.
\vskip 3mm
For the contrary, suppose  $G\ntrianglelefteq P$.
Then  by Proposition~\ref{dty1} again, $X=((A\rtimes C_1)\rtimes B)\rtimes C_2$ where $A=G_X\cong \ZZ_p^{n-1}$ and $G=A\times B$.
Set $B=\lg b\rg$. Consider the quotient group
$\ox=X/A=(\oc_1\rtimes\ob)\rtimes\oc_2$. Since we   $\oc\leq C_{\ox}(\oc_1)$ and
$\olb\notin C_{\ox}(\oc_1)$, we get  $\oc=C_{\ox}(\oc_1)\lhd X$ and so $\ob$ normalizes $\oc$. Therefore
  $B$ normalizes $A\rtimes C$, that is  $X=(A\rtimes C)\rtimes B$.

  Note that  $X=\lg \s, \ell\rg $, where $\ell$ is an involution.
  Since $\overline\ell$ is of order 2 in  $\overline X=\overline C\rtimes\overline B$ and $B\cong \ZZ_p$,
  we get  $\overline \ell=\overline\s^{\frac{kp^m}2}$. Then
$$\ox=\lg \overline{\s},\overline{\ell}\rg=\lg \overline{\s}, \overline\s^{\frac{kp^m}2}\rg=\lg \overline{\s}\rg=\oc\ne \ox,$$
  a contradiction. \qed

\begin{lem}\label{normal4}
The group $X$ contains a normal subgroup $T$ acting regularly on vertices so that $X=T\rtimes C\leq\AGL(n,p)$, where either $T=G$ or $T\cap G\cong \ZZ_p^{n-1}$ where $n\ge 2$.
\end{lem}
\demo It is easy to check that  $G\rtimes \lg z\rg =\O_1(P)$.  Set $Z=Z(G\rtimes \lg z\rg)$. It follows from
$$Z\char~\O_1(P)\char~P\char~X$$
that  $Z\char X$, which implies    $\s_2$ normalizes $Z$.
Since $C_X=1$, we know that   $[z, G]\ne 1$. Now we  have the following two cases:

\vskip 3mm
Suppose $Z\cong \ZZ_p^{s}$, where $s<n-1$. Then $G$ is the unique elementary abelian $p-$subgroup of order $p^n$ in $G\rtimes \lg z\rg$, which implies
 $$G\char G\rtimes \lg z\rg \lhd X,$$
and then  $X=T\rtimes C\le \AGL(n,p)$, where $T=G$.

\vskip 3mm Suppose $Z\cong \ZZ_p^{n-1}$. If $G\lhd X$, then $X=T\rtimes C\le\AGL(n,p)$, where $T=G$, as desired.
In what follows, we suppose $G\ntrianglelefteq X$.
Then $V_1:=\O_1(P)/Z=\lg\ola,\olz\rg\cong \ZZ_p^2$, for some $a\in G\setminus Z$. View $V_1$  as a  2-dimensional space and
consider the nontrivial projective  transformation  $\s_2$ on $\PG(V_1)$  induced by the conjugacy action.
Remind that every element in $\PGL(2,p)$ fixed at most two points.
Since $p \nmid o(\s_2)$ and $\s_2$ fixes $\lg \olz \rg$,
we know that $\s_2$ fixes exactly two points.  Therefore,
  $\s_2$  fixes another subgroup of order $p$, say  $\lg \overline {az^i}\rg $ for some  $i\in \ZZ_p^*$ and $a\ne 1$.
Then $\s_2$ fixes setwise $T:=\lg Z, a z^i\rg\cong\ZZ_p^n $.
Definitely,   $T$ is fixed by $\s_1$  too.
So $T\lhd X$. Then $X=T\rtimes C\le \AGL(n,p) $ and   $T\cap G\cong \Z_p^{n-1}$ where $n\ge 2$.
\qed

\section{Classification  Theorem}
Let $\s$, $\s_1$, $\s_2$, $C$, $C_1$ and $C_2$ be that in Section.
If $p=2$, then $X=G\rtimes C$, that is $G=T$.  Suppose that $p$ is odd. Then from Lemma~\ref{normal1},
 we know that for $X=(G\rtimes C_1)\rtimes C_2$,
 there exists a group $T\cong\ZZ_p^n$ such that $X=T\rtimes C\le \AGL(n,p)$ and $\s$ is also a skew morphism of $T$, where either $G=T$ or $G\cap T\cong \ZZ_p^{n-1}$.
 As mentioned in Section 1,  from now on  we identify $T$ with an $n$-dimensional space $V(n,p)$ and  $\s$ with a linearly transformation in $\bGL(n,p)$ so that
$\s^{-1}t\s=\s(t)$, for any $t\in T$.
 In this section, we prove shall   Theorem~\ref{main2}.
 \vskip 3mm
\f {\bf Theorem 1.2} {\it  Every regular  Cayley map  of $G\cong \ZZ_p^n$  is isomorphic to
$\MM(T\rtimes \lg \s\rg ;\, \s,\, t(-I))$,
where $\s \in {\bf M}(n,p)$ and $t$ may be given from  $T$ such that     $\{\s^i(t)\di 0\le i\le n-1\}$ generate $T$.
Moreover, given $p$ and $n$, different  $\s $ in   ${\bf M}(n,p)$ give nonisomorphic   maps.}
\vskip 3mm

 Now we have  $X=T\rtimes C$, where $C=\lg \s\rg \cong \ZZ_{kp^m}$ where $p\nmid k$ and $k$ is even if $p\ne 2$; and  also $X=\lg r,\ell\rg$, where  $o(r)= kp^m$  and  $o(\ell)=2$.
Since $\lg \s\rg$ is defined to be a point-stabilizer,   we may choose $r=\s$.   Set $\s_3=\s^{\frac{kp^m}2}$ if $o(\s)$ is even.
The proof of Theorem~\ref{main2} consists of the following three lemmas.
\begin{lem} \label{p-o1}   $\ell=t(-I)$ for some $t\in T$ such  that $T=\lg t\rg^X$; if $p$ is odd then    $\s^{\frac{kp^m}2}=-I$.
 \end{lem}
\demo (1) Suppose that $p$ is odd.  Since $X=\lg r, \ell\rg =TC$ and   $r=\s$, we may set  $\ell=t\s_3$, where  $t\in T\setminus \{1\}$. Form $\ell^2=1$ we get $\s_3(t)=-t$.
 Then  we have
 $$X=\lg r, \ell\rg =\lg \s, t\s_3\rg =\lg \s, t\rg.$$
 Set $N:=\lg t\rg ^X$,   the normal closure of $\lg t\rg$ in $X$. Then
$N=\lg \s^i(t)\di i\in \ZZ_{kp^m}\rg \leq T.$  So
$$X\le N\rtimes \lg \s\rg \le T\rtimes \lg \s \rg =X.$$
Comparing  the order of $N$ and $T$,  we get $N=T$. Now,
since $\s_3(t)=-t$ and $T=\lg \s^i(t)\di i\in \ZZ_{kp^m}\rg$, we get $\s_3(\s^i(t))=\s^i(\s_3(t))=-\s^i(t)$,  that is $\s_3=-I$.

\vskip 3mm
(2) Suppose that $p=2$.  Then $-I=I$.
If  $o(\s)$ is odd, then   $\ell=t$ for some $t\in T$ so that $X=\lg t, \s\rg$.
\vskip 3mm
Suppose $o(\s)$ is even and  $\ell\not\in T$, that is $\ell=t\s_3$ for some $t\in T$.
Then
$$X=\lg\s,t\s_3\rg=\lg\s,t\rg.$$
From  $\ell^2=1$, we get $[t,\s_3]=1$, which implies $t\in C_X(\s_3)$. Noting $\s\in C_X(\s_3)$, we get $X=\lg\s,t\rg\leq C_X(\s_3)$ and so $\lg\s_3\rg\le Z(X)$, contradicting to $C_X=1$.
Then $\ell\in T$ so that  $X=\lg\s,t\rg$.
\qed

\begin{lem} \label{p-odd} Given $X=T\rtimes \lg \s\rg $, every regular Cayley map $\MM(X; \s, \ell)$ of $\ZZ_p^n$  is isomorphic to
  $\MM(X, \s, t(-I))$, where $t$ is a  given element in $T$   such that $\lg t, \s\rg =X$.
\end{lem}
\demo  To prove the lemma we need to determine the representatives of orbits of $(\Aut(X))_\s$ on the   involutions $\ell $ such that $\lg \s, \ell \rg =X$.
\vskip 3mm
 Pick up a $t$ in $T$ such that  $\lg t, \s \rg=X$.  For any  $\ell=t_2(-I)$, where $t_2\in T$, we  assume
 $t_2=i_0t+i_1\s(t)\cdots +i_{n-1}\s^{n-1}(t)$ for some $i_0, i_1, \cdots, i_{n-1}\in\ZZ_p$.
Set  $f(\s)=i_0 I+i_1\s+\cdots +i_{n-1}\s^{n-1}.$  Then $f(\s)(t)=t_2$ and $\s f(\s)=f(\s)\s$. So
$$\begin{array}{lll} &&(t_2, \s(t_2), \cdots , \s^{n-1}(t_2))\\
&=&(f(\s)(t), \s f(\s)(t), \cdots , \s^{n-1} f(\s)(t))\\
&=&f(\s)(t, \s (t), \cdots , \s^{n-1}(t)).\end{array}$$
Since $X=\lg t,\s\rg=\lg t_2, \s\rg$ and both
$\{t, \s (t),\cdots,\s^{n-1}(t)\}$ and $\{t_2, \s (t_2),\cdots,\s^{n-1}(t_2)\}$ are bases of $T$,
we  know that $f(\s)\in \bGL(n,p)$, which commutes with $\s$.  So, $f(\s)$  induces an automorphism on $X$ which  fixes $\s$ and maps $t$ to $t_2$.
Therefore, we may choose $\ell=t(-I)$ for this given $t$.
\qed

\vskip 3mm

\begin{lem}\label{def1} Given $\ZZ_p^n$, different $\s$  in   ${\bf M}(n,p)$ give nonisomorphic   maps.
\end{lem}
\demo  Suppose that $$\MM(T\rtimes \lg \s_1\rg ;\, \s_1,\, t_1(-I))\cong  \MM(T\rtimes \lg \s_2 \rg ;\, \s_2,\, t_2(-I)),$$
for $\s_1, \s_2$ in ${\bf M}(n,p)$. Then  there exists an isomorphism $\psi: T\rtimes \lg \s_1\rg \to  T\rtimes \lg \s_2\rg $ such that
$\psi(\s_1)=\s_2$ and $\psi(t_1(-I))=t_2(-I)$. Then
$$\psi(-I)=\psi(\s_1^{\frac{kp^m}2})=\s_2^{\frac{kp^m}2}=-I,$$
which implies $\psi(t_1)=t_2$.
Therefore
$$\psi(T)=\psi(\lg \s_1^i(t_1)\di i\in \ZZ_{kp^m}\rg)=\lg \s_2^i(t_2)\di i\in \ZZ_{kp^m}\rg =T.$$
By Proposition~\ref{dt1}, $\psi$ can be induced by  an inner automorphism  $\Inn(g)$ of $\AGL(n,p)$,  for some $g\in \bGL(n,p)$, in other words,
  $\s_1$ and $\s_2$ are conjugate in $\bGL(n,p)$, which from the definition of ${\bf M}(n,p)$  implies $\s_1=\s_2$.
  Therefore,  different $\s$  in ${\bf M}(n,p)$ give nonisomorphic maps.\qed

\section{${\bf M}(n,p)$ and ${\bf M}_1(n,p)$}
From Theorem~\ref{main2}, we know that the  regular Cayley maps of $\ZZ_p^n$ are determined by   ${\bf M}(n,p)$.  So in this section,  we shall  determine ${\bf M}(n,p)$ and ${\bf M}_1(n,p)$ as well,
 that is  Theorem~\ref{main3}. The proof for that  consists of Lemmas~\ref{Mnp} and ~\ref{unbla}.     Some technics in polynomial theory over $\bF_p$ will be used.
 \vskip 3mm
 For any $\s\in {\bf M}(n,p)$, let $$f(x)=x^n-r_{n-1}x^{n-1}-\cdots-r_1x-r_0$$
 be the minimal polynomial of $\s$. In fact, it is the characteristic polynomial of $\s$.

\vskip 3mm
 Write
 \begin{eqnarray}\label{f}
 f(x)=\prod_{i=0}^{s-1}p_i(x)^{k_i},
 \end{eqnarray}
where every $p_i(x)$ is monic and irreducible on $\bF_p$. Then $$\s =\s_0\s_1\cdots \s_{s-1}\quad {\rm and}\quad   T=T_0\bigoplus T_1\bigoplus\cdots \bigoplus T_{s-1},$$
where  the minimal polynomial of $\s_i$ on $T_i$ is precisely $p_i(x)^{k_i}$ and $\s_i\di_{T_j}=1$  for any $j\ne i$.
\vskip 3mm
Let $\a_i$ be  a root of $p_i(x)$, where  $0\le i\le s-1$ and $u_i$ a order of $\a_i$ in $\bF_p(\a)^*$.
Let  $t_i$  be the smallest integer with $p^{t_i}\geq k_i$.    Remind that for any $e\in\ZZ^+$ and a prime $p$,  by $[e]_p=s$, we mean that   $p^s\di e$ but $p^{s+1}\nmid e$.

Then we have
\begin{lem}\label{2,p}  Let $f(x)$ be as in  Eq(\ref{f}). Then
 \begin{enumerate}
   \item[\rm(1)] $o(\s_i)=u_ip^{t_i}$ for any $i$ and $o(\s)=\ord(f(x))=\lcm\{ u_ip^{t_i}\di 0\le i\le s-1\}$;
   \item[\rm(2)]   Suppose $p$ is odd. Then $\lg \s \rg $ contains $-I$ if and only if   $f(x)|(x^{\frac{o(\s)}2}+1)$, and if and only if
   $p_i(x)^{k_i}\di (x^{\frac{o(\s_i)}2}+1)$    and
   $[u_i]_2$ is a constant for all $0\le i\le s-1$.
   \end{enumerate}
\end{lem}
\demo
(1)  Recall that the order $\ord(p_i(x))$ of $p_i(x)$ is defined to be the minimal $e_i$
such that $p_i(x)\di (x^{e_i}-1)$. By Proposition~\ref{RN2},  $\ord(p_i(x))=u_i:=\o(\a_i)$ in $\bF_p(\a)^*$, and
 $\ord(p_i(x)^{k_i})=u_ip^{t_i}.$
Since $p_i(x)\
^{k_i}$ is the minimal polynomial of $\s_i$ and $u_ip^{t_i}$ is the minimal number
such that $p_i(x)^{k_i}\di (x^{u_ip^{t_i}}-1)$, $u_ip^{t_i}$ is the minimal number
such that $\s_i^{u_ip^{t_i}}=1$, which implies $o(\s_i)=u_ip^{t_i}=\ord(p_i(x)^{k_i})$ for any $i$.
Since $\lg \s_i\rg \cap \lg \s_j\rg =1$ and $[\s_i, \s_j]=1$, for any $i\ne j$, we know from Proposition~\ref{RN2}  that
$$o(\s)=o(\Pi_{i=0}^{s-1}\s_i)=\lcm\{o(\s_i)\di 0\le i\le s-1\}=\lcm\{ u_ip^{t_i}\di 0\le i\le s-1\}=
\ord(f(x)).$$

(2)  Suppose $p$ is odd. Set $e=o(\s)=\ord(f(x)).$  Write  $\ord(p_i(x)^{k_i})=u_ip^{t_i}=2^{m_i}h_i$, where $h_i$ is odd.  Without loss of any generality,  assume that  $m_0\geq m_1\geq\cdots\geq m_{s-1}$.
    Set  $h=\lcm(h_0,h_1, \cdots, h_{s-1})$. Then  by Proposition~\ref{RN2}, we get $$e=\lcm(2^{m_0}h_1, 2^{m_1}h_2, \cdots, 2^{m_{s-1}}h_{s-1})=2^{m_0}h.$$
\vskip 3mm
Now,   $-I\in\lg\s\rg$ if and only if  $e$ is even and  $\s^{\frac e2}+I=0$, if and only if
   $f(x)|(x^{(\frac{o(\s)}2)}+1)$, which is the first part of (2).

\vskip 3mm
Suppose that $f(x)\mid (x^{\frac e 2}+1)$.
Then $[e]_2\geq1$ and $p_i(x)^{k_i}\mid (x^{{\frac e 2}}+1)$. First we show that $[u_i]_2$ is a constant for all $0\le i\le s-1$.
For the contrary, suppose $m_{j}\lneqq m_{j+1}$ for some $j\ge 1$. Since $2^{m_{j}}h_{j}\di \frac e2$,
we have $(x^{2^{m_j}h_j}-1)|(x^{\frac e2}-1)$, which implies $p_{j}(x)^{k_{j}}|(x^{\frac e2}-1)$.
Since $(x^{\frac e 2}-1,x^{\frac e 2}+1)=1$, we get   $p_{j}(x)^{k_{j}}\nmid (x^{\frac e2}+1)$, which implies  $f(x)\nmid(x^{\frac{e}{2}}+1)$, a contradiction.
Therefore, $[u_i]_2=m_i=[e]_2$ for any $0\leq i\leq s-1$.
Since  $$p_i(x)^{k_i}\mid (x^{2^{[e]_2}h_i}-1), p_i(x)^{k_i}\mid (x^{{\frac e 2}}+1),
(x^{2^{[e]_2}h_i}-1,x^{{\frac e 2}}+1)=x^{2^{[e]_2-1}h_i}+1=x^{\frac{o(\s_i)}2}+1,$$ we get
$$p_i(x)^{k_i}\mid (x^{\frac{o(\s_i)}2}+1).$$

\vskip 3mm
Conversely, suppose that $p_i(x)^{k_i}\mid (x^{\frac{o(\s_i)}2}+1)$ and $s_i=[e]_2$ for any $0\leq i\leq s-1$.
Since  $(x^{2^{m_i-1}h_i}+1)\mid (x^{\frac e2}+1)$ as $m_i=[e]_2$, we get  $p_i(x)^{k_i}\mid (x^{\frac e2}+1)$. Therefore, $f(x)\mid (x^{\frac e2}+1)$.\qed

\begin{lem}\label{2,p,1}  Let $f(x)$ be as in  Eq(\ref{f}), where $p$ is odd and $p\di o(\s)$. Then
    $\lg \s^{\frac{o(\s)}p} \rg $  fixes pointwise
   an $(n-1)$-dimensional subspace if and only if
   there exists an $i$ such that $p_i(x)=x-\a_i$
   such that $t_i\gneqq t_j$ for any $j\ne i$ and
   $k_i=p^{t_i-1}+1$.
\end{lem}
\demo (1)  First we consider a special case:  $f(x)=p(x)^k$ is the minimal polynomial of $\s$, where $p(x)$ is an irreducible polynomial and $\ord(f(x))=\ord(p(x))p^{t}$, where $t\ge 1$.

\vskip 3mm

Suppose $\lg \s^{\frac{o(\s)}p} \rg $ fixes pointwise an $(n-1)$-dimensional subspace.
Then   $\s$ fixes a 1-dimensional subspace $(\s^{\frac{o(\s)}p}-I)(T)$,
 which implies that $\s$ has an eigenvalue, say $\a$ in $\bF_p$. Since  $p(x)$ is irreducible, we get $p(x)=x-\a$ so that $f(x)=(x-\a)^k.$
 Set   $o(\a)=u$ so that $o(\s)=o(\a)p^t$ and  $\s^{\frac{o(\s)}p}=\s^{up^{t-1}}$.
 Set $\s=\a I+\t $. Then $\t^k=(\s-\a I)^k=0$, which means that  $\t$ is nilpotent. Then
 $$\s^{up^{t-1}}=(\a I+\t )^{up^{t-1}}=(\a I+\t^{p^{t-1}})^u=I+u\a^{u-1}\t^{p^{t-1}}+\cdots +\t^{up^{t-1}}$$
So  $$(\s^{up^{t-1}}-I)T=(u\a^{u-1}\t^{p^{t-1}})T+\cdots +(\t^{up^{t-1}})T=\t^{p^{t-1}}T,$$
which is 1-dimensional if and only if $\t(\t^{p^{t-1}})T=0$.  By the definition of $k$ and $p^t$, we get $1+p^{t-1}=k$.

\vskip 3mm
(2) Coming back to the general case,
suppose $p\di o(\s)$ and $\lg \s^{\frac{o(\s)}p} \rg $ fixes pointwise a $(n-1)$-dimensional subspace.
Recall  $f(x)=\Pi_{i=0}^{s-1}p_i(x)^{k_i},$
where every $p_i(x)$ is monic and irreducible on $\bF_p$ and $\ord(p_i(x)^{u_i})=\ord(p_i(x))p^{t_i}$,
and
$$\s =\s_0\s_1\cdots \s_{s-1}\quad {\rm and}\quad   T=T_0\bigoplus T_1\bigoplus\cdots \bigoplus T_{s-1},$$
where $\s_i\di_{T_j}=1$ for any $j\ne i$
and the minimal polynomial of $\s_i$ on $T_i$ is precisely $p_i(x)^{k_i}$
and $\deg(p_i(x)^{k_i})=\dim(T_i)$.
Without loss of any generality,  suppose that $t_0\geq t_1\geq\cdots\geq t_{s-1}$.
Now we are proving   $t_0\gneqq t_1$. For the contrary,  suppose that $t_0=t_1$. Then both  $\s_0^{\frac{o(\s)}p}$ and  $\s_1^{\frac{o(\s)}p}$  are nontrivial elements in
$\lg\s^{\frac{o(\s)}p}\rg$,
which implies that for $i=0$ and 1,  $\s_i^{\frac{o(\s)}p}$ fixes pointwise a
 subspace in $T_i$ of dimension at most $(\dim(T_i)-1)$. Therefore $\s^{\frac{o(\s)}p}=\s_0^{\frac{o(\s)}p}\s_1^{\frac{o(\s)}p}\cdots \s_{s-1}^{\frac{o(\s)}p}$
  fixes a pointwise subspace  of dimension at most $(n-2)$, a contradiction. Thus we get
   $\s^{\frac{o(\s)}p}=\s_0^{\frac{o(\s)}p}=(\s_0^{\frac{o(\s_0)}p})^{\frac{o(\s)}{o(\s_0)}}$,
where $p\nmid \frac{o(\s)}{o(\s_0)}$. Pick up $j$ such that $\frac{o(\s)}{o(\s_0)}j\equiv 1(\mod p)$. Then
 $\s_0^{\frac{o(\s_0)}p}=\s^{{\frac{o(\s)}p}j}$ fixes pointwise an $(n-1)$-dimensional subspace.
Form (1),  we get $p_0(x)=x-\a_0$ for some $\a_0$ and $k_0=p^{t_0-1}+1$, as desired.
\vskip 3mm
Conversely, suppose that $p_0(x)=x-\a_0$ such that $t_0\gneqq t_j$ for any $j\ne 0$ and $k_0=p^{t_0-1}+1$.
Then $\s^{\frac{o(\s)}p}=\s_0^{\frac{o(\s)}p}$ and  by (1),
$\lg \s^{\frac{o(\s)}p}\rg=\lg \s_0^{\frac{o(\s)}p}\rg=\lg \s_0^{\frac{o(\s_0)}p}\rg$, which
fixes pointwise an $(n-1)$-dimensional subspace.
\qed
\vskip 3mm
Let $p$ be odd. Recall that ${\bf MF}_p$ denotes the set of monic polynomials  over $\bF_p$.
In observing of Lemma~\ref{2,p},
to characterize  all the polynomials of even order we need to introduce the following set:
for any $i, j\in \ZZ^+$ and $2^j\di (p^i-1)$, let
\begin{eqnarray} \label{Aij}
A_{ij}=\{f(x)|\,f(x)\in {\bf MF}_p[x], \deg(f(x))=i,f(x)~\rm{is~irreducible},\, [\ord(f(x))]_2=j\},
\end{eqnarray}
where $1\le j\le [p^i-1]_2$. Note that $A_{ij}=\emptyset ,$ for some $p, n, i$ and  $j$. The value $|A_{ij}|$ will be computed in next section and here  one of its properties  is given.
\begin{lem}\label{prop-of-Aij}
For any  $f(x)\in A_{ij}$, we have $f(x)|(x^{\frac 12\ord(f(x))}+1)$ and
$f(x)^h|(x^{\frac 12 \ord(f(x)^h)}+1)$, where $h\in\ZZ^+$.
\end{lem}
\demo
Set $e=\ord(f(x))$. By Proposition~\ref{RN2},  $\ord(f(x)^{h})=ep^t$, where $t$ is the smallest positive integer such that $p^t\geq h.$
Since $j\ge 1$, we know $e$ is even. Then $f(x)\mid (x^{\frac e2}+1)(x^{\frac e2}-1)$. Note that
$f(x)\nmid x^{\frac e2}-1$ as $\ord(f(x))=e$. Hence $f(x)\mid x^{\frac e2}+1$, as $f(x)$ is irreducible.
Since $f(x)^h|(x^{ep^t}-1)$ and $(x^{\frac {ep^t}2}+1,x^{\frac {ep^t}2}-1)=1$, we get
 $f(x)^h|(x^{\frac{ep^t}{2}}+1)$.
\qed
\vskip 3mm
Since there is a one-to-one correspondence between $\s$ and its  monic  minimal polynomial, we introduce
$$M(n,p)=\{f(x)\in {\bf MF}_p[x]\di f(x)\, {\rm is\, the\,   minimal\, polynomial\, of}\, \s \, {\rm in}\,  {\bf M}(n,p)\},$$
$$M_1(n,p)=\{f(x)\in {\bf MF}_p[x]\di f(x)\, {\rm is\, the\,   minimal\, polynomial\, of}\, \s\, {\rm in}\, {\bf M}_1(n,p)\}.$$
 Let $d=\max\{[p^i-1]_2\di 1\le i\le n\}$. If $2^j\nmid (p^i-1)$,  then  we   define  $A_{ij}=\emptyset $.
 Then for any given $1\le i\le n$ and $1\le j\le d$,  $A_{ij}$ and $|A_{uj}|$  are uniquely  defined.
 Then  $M(n,p)$  and $M_1(n,p)$ may be divided into  subsets   which are characterized by  $j=[\ord(f(x)]_2.$

 \vskip 3mm
 Now we are ready to prove Theorem~\ref{main3}, by combining Lemmas~\ref{Mnp}  and ~\ref{unbla}.
\begin{lem}   \label{Mnp} Let $p$ be odd. Then we have
\begin{eqnarray}\label{mnp}
 M(n,p)=\bigcup_{j=1}^d \{\prod_{i=1}^{n}\prod_{r=1}^{|A_{ij}|}f_{ijr}(x)^{k_{ijr}}\mid f_{ijr}(x)\in A_{ij}\ne \emptyset,
\sum_{i=1}^{n}(\sum_{r=1}^{|A_{ij}|}k_{ijr})i=n\},
\end{eqnarray}
\begin{eqnarray}\label{mnp1}
M_1(n,p)=\{f(x)\in M(n,p)\mid  \,{\rm there\, exists\, an\, unique\, } (j,r)\, {\rm  such\, that\, } \end{eqnarray}
 $$k_{1jr}=p^{[o(\s)]_p-1}+1, \, {\rm but}\,  [k_{ijr}]_p\lneqq [o(\s)]_p, \,{\rm for\, other\, } k_{ijr}\}.$$

 \end{lem}
\demo  (1)\, Since $2^j\di (p^i-1)$, we know that $1\le j\le [p^i-1]_2$ and so $1\le j\le d$.
 It follows from  $\sum_{i=1}^{n}(\sum_{r=1}^{|A_{ij}|}k_{ijr})i=n$ that  for any $j$ such that $A_{ij}\ne \emptyset$,  every polynomial  $\prod_{i=1}^{n}\prod_{r=1}^{|A_{ij}|}f_{ijr}(x)^{k_{ijr}}$  in right side  of Eq(\ref{mnp}) is of degree $n$.
 Since  $f_{ijr}(x)\in A_{ij}\ne \emptyset$, it follows from Lemma~\ref{prop-of-Aij} that $f_{ijr}(x)|(x^{\frac 12\ord(f_{ijr}(x))}+1)$ and
$f_{ijr}(x)^{k_{ijr}}|(x^{\frac 12 \ord(f_{ijr}(x)^{k_{ijr}})}+1)$. Moreover,  $[\ord (f_{ijr}(x)^{k_{ijr}}]_2=j$ for any $1\le r\le |A_{ij}|$.  By Proposition~\ref{2,p}.(2),
 $-I\in \lg \s \rg $. Therefore, $\prod_{i=1}^{n}\prod_{r=1}^{|A_{ij}|}f_{ijr}(x)^{k_{ijr}}\in M(n,p).$

 \vskip 3mm
 Now we are proving that every polynomial  $f(x)$ in $M(n,p)$ has the form of the right side of  Eq(\ref{mnp}).   Set $e=o(\s)=\ord(f(x)).$  By Lemma~\ref{2,p}, write  $\ord(p_i(x)^{k_i})=u_ip^{t_i}=2^{m}h_i$, where $m\ge 1$.
    Write  $f(x)=\Pi_{i=0}^{s-1}p_i(x)^{k_i}$ and set  $h=\lcm(h_0,h_1, \cdots, h_{s-1})$. Then  by Lemma~\ref{2,p}, we get $e=2^mh.$
 Then $p_i(x)\in A_{\deg(p_i(x)), m}$ and $k_i=k_{\deg(p_i(x)), m, r}$ for some $1\le r\le |A_{\deg(p_i(x)), m}|$.
 Moreover,  since  $2^m\di (p^{\deg(p_i(x))}-1),$  we have $1\le m\le [p^{\deg(p_i(x))}-1]_2$.
  Enumerating  the  dimension $n$ and relabeling $\deg(p_i(x))$ and $m$ by $i$ and $j$, respectively,  we get
$$\sum_{i=1}^{n}(\sum_{r=1}^{|A_{ij}|}k_{ijr})i=n.$$

\vskip 3mm
(2)\, By the definition, a $\s$ in  ${\bf M}(n,p)$ is contained in ${\bf M}_1(n,p)$ if and only if  $p\di o(\s)$ and $\lg \s^{\frac{o(\s)}p} \rg $ fixes pointwise
   an $(n-1)$-dimensional subspace. By Lemma~\ref{2,p,1}, the later holds if any only if  the minimal polynomial $f(x)$ of $\s$  contains a factor   $p_i(x)=x-\a_i$
   such that  $k_i=p^{t_i-1}+1$ and $t_i\gneqq t_j$ for any $j\ne i$. In other words, $[o(\s)]_p=t_i$.  Therefore, $M_1(n,p)$ is given by Eq(\ref{mnp1}). \qed

\begin{lem} \label{unbla}
 The case  $T\cap G\cong \Z_p^{n-1}$ in Theorem~\ref{main1} happens if and only if
 $p\ne 2$ and $\s \in {\bf M}_1(n,p)\ne \emptyset,$ equivalently, $f(x)\in M_1(n,p)$, \
 where $f(x)$ is the minimal polynomial of $\s$.
 \end{lem}
 \demo Suppose $T\cap G\cong \Z_p^{n-1}$. Then  from the proof of Lemma~\ref{normal4}, we know that  $p\ne 2$ and $Z(T\rtimes\lg\s^{\frac{o(\s)}p}\rg)\cong\ZZ_p^{n-1},$ equivalently,
 $\lg\s^{\frac{o(\s)}p}\rg$  fixes pointwise an $(n-1)$-dimensional subspace in $T$, which implies $\s \in {\bf M}_1(n,p)\ne \emptyset,$ that is,  the minimal polynomial $f(x)$ of $\s$
 is contained in $M_1(n,p)$.

\vskip 3mm
Conversely, for any $f(x)\in M_1(n,p)$, we get $f(x)^{\frac{\ord(f(x))}p}=(x-\a)^{p^t+1}$. In other words, $\s^{\frac{o(\s)}p}$    fixes pointwise an  $(n-1)$-dimensional subspace in $T$. Therefore,
 $Z(T\rtimes\lg\s^{\frac{o(\s)}p}\rg)\cong\ZZ_p^{n-1},$ which implies that   $T\cap G\cong \Z_p^{n-1}$.\qed

\section{Enumeration   Theorem}
By Lemma~\ref{Mnp},  $M(n,p)=\bigcup_{j=1}^dM_j$, where
\begin{eqnarray} \label{mj} M_j=\{\prod_{i=1}^{n}\prod_{r=1}^{|A_{ij}|}f_{ijr}(x)^{k_{ijr}}\mid f_{ijr}(x)\in A_{ij}\ne \emptyset,
\sum_{i=1}^{n}(\sum_{r=1}^{|A_{ij}|}k_{ijr})i=n\}. \end{eqnarray}
So,  to enumerate $|{\bf M}(n,p)|=|M(n,p)|$, we need to know  $|A_{ij}|$.

\begin{lem}\label{number-of-Aij} Given $i, j\in \ZZ^+$, set
$$I_{ij}=\{e|[e]_2=j,e|(p^i-1),\, {\rm but}\, e\nmid(p^{i'}-1),~\forall~{\rm proper\, divisor }\, i'\, {\rm of}\, i \}.$$
 Then  \begin{eqnarray} \label{Aij1}
 |A_{ij}|=\frac 1i \sum_{e\in I_{ij}}\varphi(e).
 \end{eqnarray}
\end{lem}
\demo Set   $$L_{ij}=\{\ord(f(x))|f(x)\in A_{ij} \}.$$
By Proposition~\ref{RN3}, the number of monic irreducible polynomials in $\bF_p[x]$
of degree $i$ and order $e$ is equal to $\frac{\varphi(e)}i$.
Then $|A_{ij}|=\frac 1i \Sigma_{e\in L_{ij}}\varphi(e)$.
Note that  for any $f(x)\in A_{ij}$, we have $\ord(f(x))=2^jk$ for some odd integer $k$.
Then  we only need  to show
$L_{ij}=I_{ij},$  as defined in the lemma.

\vskip 3mm
For any $f(x)\in A_{ij}$, by Proposition~\ref{RN4}, $\ord(f(x))$ is equal to the order of any root of $f(x)$ in the multiplicative group $\bF_{p^i}$. Hence
$$\begin{array}{lll}
L_{ij}&=&\{\ord(f(x))|f(x)\in A_{ij} \}\\
      &=&\{e|o(\b)=e,\,\b~\rm{is~the~root~of}~\textit{f(x)},\,\, \textit{f(x)}\in A_{ij}\}\\
      &=&\{e|o(\b)=e,\,[e]_2=j,\,\bF_p(\b)=\bF_{p^i}\}.
  \end{array}$$
If $\bF_p(\b)=\bF_{p^i}$, then $\b$ is contained in $\bF_{p^i}$ but not in any proper subfield of $\bF_{p^i}$, which implies
$o(\b)|(q^i-1)$ and $o(\b)\nmid(q^{i'}-1)$, for any proper divisor $i'$ of $i$.
Therefore, we get
$$L_{ij}=\{e|[e]_2=j,e|(q^i-1),\, {\rm but}\, e\nmid(q^{i'}-1),~\forall~{\rm proper\, divisor }\, i'\, {\rm of}\, i \}=I_{ij}.$$
\qed

\vskip 3mm
Now  we are ready to enumerate  $|M(n,p)|$. Recall $d=\max\{ [p^i-1]_2\di 1\le i\le n\}$.

\begin{lem} \label{nMnp}
Given $1\le j\le d $,
set $a_{ij}:=|A_{ij}|$ and
let  $n_j$ be the number of the solutions $\{k_{ijs}\di 1\le s\le a_{ij}\}$  of
the following integer equation over $\ZZ^+$ such that $a_{ij}\ne 0$:
\begin{eqnarray} \label{n5}
 \sum_{i=1}^{n}(\sum_{r=1}^{a_{ij}}k_{ijr})i=n.
\end{eqnarray}
Then $|M(n,p)|=\sum_{j=1}^d n_j$.
\end{lem}
\demo   Set $a_{ij}:=|A_{ij}|$ and write
$$A_{ij}=\{f_{ij1}(x), f_{ij2}(x), \cdots, f_{ij a_{ij}}(x)\}.$$
Then every polynomial  $f(x)$ in $M_j$ (defined in Eq(\ref{mj}))  can  be uniquely decomposed as
$$f(x)=\prod_{i=1}^{n}\prod_{r=1}^{a_{ij}}f_{ijr}^{k_{ijr}},$$
where some $k_{ijr}$ may be zero. Hence, this $f(x)$   is uniquely determined by the parameters  $$\{k_{ijr}\di 1\le i\le n,\,   1\le r\le a_{ij}\},$$
which is a solution of Eq(\ref{n5}). Conversely, every solution of  Eq(\ref{n5}) gives a polynomial $f(x)$, which is contained   in $M_j$.
 Therefore, the cardinality of $M_j$ is $n_j$ and so
 $|M(n,p)|=\sum_{j=1}^d n_j$, as desired.
\qed

\vskip 3mm
The following  example makes an illustration for   Lemma~\ref{nMnp}.

\begin{exam}\label{2,3,4}
 $|{\bf M}(2,3)|=4$  and   $|{\bf M}(3,3)|=5$.
\end{exam}
\demo (1)\,  $n=2$ and $p=3$:  By the definition of $i$ and $j$, we have $2^j\di (3^i-1)$ where $1\le i\le 2$ and $j\ge 1$ for which there exist possible nonempty    $A_{ij}$,  so $$(i,j)\in \{ (1, 1), (2, 1), (2, 2), (2, 3)\}.$$
Check that
$$A_{11}=\{ x+1\};~A_{21}=\emptyset;~A_{22}=\{x^2+1\};~A_{23}=\{x^2+x-1,x^2-x-1\}.$$
So
$$a_{ij}:=|A_{ij}|=\left\{\small\begin{array}{ll}
    1,  & i=1,\,j=1;  \\
    1,  & i=2,\,j=2;  \\
    2,  & i=2,\,j=3;  \\
    0,  & \rm{others}.\\
    \end{array}\right..$$
    Suppose that $j=1$: Eq(\ref{n5}) is :
    $$k_{111}\cdot 1=2,$$
    which has  only one solution, that is $n_1=1$.   Suppose that $j=2$: solving the equation:
    $$2k_{221}=2,$$
    we get only one solution and so $n_2=1$. Suppose that $j=3$: solving the equation:
    $$2(k_{231}+k_{232})=2,$$
    we get  two solutions and so $n_3=2$.
    Then by Lemma~\ref{nMnp}, we get $$|{\bf N}(2,3)|=n_1+n_2+n_3=4.$$.

\vskip 3mm (2)\, $n=3$ and $p=3$:
First, we need to calculate $a_{ij}$. One can check $a_{ij}$ for all the possibilities for $i, j$:
$$a_{ij}=\left\{\small\begin{array}{ll}
    1,  & i=1,\,j=1;  \\
    1,  & i=2,\,j=2;  \\
    2,  & i=2,\,j=3; \\
    4,  & i=3,\,j=1;  \\
    0,  & \rm{others}.\\
    \end{array}\right.$$
Then, with the same way as that in (1), we get $n_1=5,$~and $n_2=n_3=n_4=0$.
So $|{\bf N}(3,3)|=5$.\qed
\section{Genera  of the maps}
It is well-known that the genus of an orientably-regular map is given by
$$g = 1 +\frac12(|E|-|V |-|F|),$$ where   by $|V |, |E|$ and $|F|,$ we denote the number of vertices,
edges and faces of the given map (embedding), respectively.
Now, for our  regular Cayley map $\MM(X,r,\ell)$ of $\ZZ_p^n$, we have that  $|V|=p^n$ and $|E|=\frac{|X|}2$.  In what follows we shall  calculate $|F|=\frac{|X|}{|\lg r\ell\rg|}$.
By Theorem~\ref{main2}, we get $r=\s$ and $\ell$ is $t(-I)$, where $t$ is given element $T$ such that   $\lg t,\s\rg=X$,
and note $-I=I$ when $p=2$.
First  we have the following two lemmas:

\begin{lem} \label{g-p-odd}
For the regular Cayley map $\MM(T\rtimes \lg \s\rg;  \s, t(-I))$ of $\ZZ_p^n$, where $p$ is odd, set $r=\s$ and $\ell=t(-I)$. Then
$$o(r\ell)=\left\{\small\begin{array}{ll}
    \a, & {\rm if}\, k\equiv 0(\mod 4); \\
    \frac \a2,   &{\rm if}\,  k\equiv 2(\mod 4),  (x-1)^{p^m}\nmid f_1(x);\\
    \frac{p\a}2, & {\rm if}\, k\equiv 2(\mod 4), (x-1)^{p^m}\| f_1(x),\\
    \end{array}\right. $$
    where $\a=o(\s)=kp^m$, $(p,k)=1$  and $f_1(x)$ is the minimal polynomial of $\s^2$.
\end{lem}
\demo  Suppose $p$ is odd and $\a=o(\s)=kp^m$ where $p\nmid k$. Since $r\ell=\s^{1+\frac \a 2}t^{-1}$, where  $\gcd(1+\frac \a2, \a)=1$ or $2$, when $k\equiv 0(\mod 4)$ or  $k\equiv 2(\mod 4)$, respectively,   we  deal with two cases, separately:

\vskip 3mm
{\it Case 1: $k\equiv 0(\mod 4)$}
\vskip 3mm
Since $(1+\frac \a2, \a)=1$ and $T\lhd X$,  $o(r\ell)\ge \a$.  Since $$(r\ell)^{\frac{\a}2}=(\s^{1+\frac{\a}2}t^{-1})^{\frac \a 2}=(-I)t'$$
for some $t'\in T$, we get $(r\ell)^{\a}=1,$ which implies  $o(r\ell)=\a$.

\vskip 3mm
{\it Case 2: $k\equiv 2(\mod 4)$}
\vskip 3mm
 Suppose $k\equiv 2(\mod 4)$.
 Then $(1+\frac \a 2, \a)=2$ and
$$ T=\lg t\rg^{X}=\lg t\rg^{\lg \s\rg}=\lg t\rg^{(\lg-I\rg\lg \s^2\rg)}=\lg t\rg^{\lg \s^2\rg}.$$
Let $f_1(x)$ be  the minimal polynomial of $\s^2$.
Since $o(\s^2)=\frac\a2$,  $f_1(x)|(x^{\frac\a2}-1).$
Set $\s_1=\s^{1+\frac \a2}$. Then
$$(r\ell)^{\frac \a2}=(\s^{1+\frac{\a}2}t^{-1})^{\frac\a2}=(I+\s_1+\cdots +\s_1^{(\frac \a2-1)})(t)=(I+\s^2+\cdots +(\s^2)^{(\frac \a2-1)})(t),$$
which implies  that  $o(r\ell)=\frac \a2 $ if and only if $(I+\s^2+\cdots +(\s^2)^{(\frac \a2-1)})(t)$ is zero vector.
Let $$g(x)=1+x+x^2+\cdots+x^{\frac \a2-1}.$$
Since $g(\s^2)\s=\s g(\s^2)$ and $\lg t\rg^{\lg\s\rg}=T$,
it follows that $-g(\s^2)(t)=1$ if and only if $T^{\lg-g(\s^2)\rg}=1$, which implies $g(\s^2)=0$.
Hence, $o(r\ell)=\frac \a2$ if and only if $g(\s^2)=0$, and if and only if $f_1(x)|g(x)$.
Since  $(x-1)^{p^m}\parallel x^{\frac \a2}-1$ and $x^{\frac \a2}-1=(x-1)g(x)$,
we get that $o(r\ell)=\frac \a2$ or $\frac {p\a}2$  if and only if
$(x-1)^{p^m}\nmid f_1(x)$ or $(x-1)^{p^m}\| f_1(x)$, respectively.
\qed

\begin{lem} \label{g-p-2}
For the regular Cayley map $\MM(T\rtimes \lg\s\rg;\, \s, t)$ of $\ZZ_2^n$, set $r=\s$ and $\ell=t$.
Then
$$o(r\ell)=\left\{\small\begin{array}{cc}
    \a, & {\rm if}\, (x-1)^{2^m}\nmid f(x) \\
    2\a, & {\rm if}\, (x-1)^{2^m}\| f(x)\\
    \end{array}\right.,$$
where  $\a=o(\s)=k2^m$, $k$ is odd and $f(x)$ is the minimal polynomial of $\s$.
\end{lem}

\demo
Note that
$$(r\ell)^{\a}=(\s t)^{\a}=(I+\s+\cdots +\s^{\a-1})(t).$$
Then $o(r\ell)=\a$ if and only if $(I+\s+\cdots +\s^{\a-1})(t)=1$.
Let  $g(x)=1+x+\cdots+x^{\a-1}$.
Since $g(\s)\s=\s g(\s)$ and $\lg t\rg^{\lg\s\rg}=T$,
then $g(\s)(t)=1$ if and only if $g(\s)=0$. Since $f(x)$ is the minimal polynomial of $\s$,
$o(r\ell)=\a$ if and only if $f(x)\mid g(x)$. Note that $(x-1)^{2^m}\parallel x^{2^mk}-1$ and $x^{2^mk}-1=(x-1)g(x)$, we get that
$o(r\ell)=\a$ or $2\a$ if and only if $(x-1)^{2^m}\nmid f(x)$ or
 $(x-1)^{2^m}\| f(x)$, respectively.
\qed

\vskip 4mm
 With Lemmas~\ref{g-p-odd} and \ref{g-p-2},  inserting  $|V|=p^n$ and $|E|=\frac{|X|}2$ and  $|F|=\frac{|X|}{|\lg r\ell\rg|}$ in $g = 1 +\frac 12 (|E|-|V |-|F|),$  we immediately get the following lemma:

\begin{lem} \label{genus}
Let $g$ be the genus of regular Cayley map $\MM(T\rtimes \lg \s \rg ;\, \s,\, t(-I))$, where
$\s\in {\bf M}(n,p)$. Then
 $g=1+\frac14p^n(\a-2\b-2)$, where $\a=o(\s)=kp^m$, $(k,p)=1$ and
$$\b=\left\{\small\begin{array}{ll}
    1,         & {\rm if } \,(x-1)^{2^m}\nmid f(x),\, p=2;\\
    \frac 12,  &{\rm if } \, (x-1)^{2^m}\| f(x),\, p=2; \\
    1,         & {\rm if }\, k\equiv 0(\mod 4),\, p\neq2;\\
     2,        & {\rm if } \,(x-1)^{p^m}\nmid f_1(x),  k\equiv 2(\mod 4),\, p\neq2; \\
    \frac 2p,  &{\rm if } \, (x-1)^{p^m}\| f_1(x), k\equiv 2(\mod 4),\,  p\neq2,\\
    \end{array}\right.,$$
where $f(x)$ and $f_1(x)$ are minimal polynomials of $\s$ and $\s^2$, respectively.
\end{lem}

\begin{exam}\label{2,2,g}
Using Theorem~\ref{main4}, we can obtain the genera  of regular Cayley maps
$\MM(T\rtimes\lg\s\rg ;\, \s,\, t(-I))$, when $\s\in \{{\bf M}(2,2),\, {\bf M}(2,3),\,{\bf M}(3,2), {\bf M}(3,3)\}$.
\vskip 3mm

(1)\, $\s\in {\bf M}(2,2)$: {\footnotesize
$$\begin{tabular}{|c|c|c|}
  \hline
  $f(x)$ & $x^2+x+1$ & $x^2+1$\\
  \hline
  $\a$    &    $3$    &  $2$ \\
    \hline
  $\b$    &    $1$    &  $\frac12$ \\
    \hline
  $g$    &    $0$    &  $0$ \\
  \hline
\end{tabular}
$$}
\vskip 3mm
(2)\, $\s\in {\bf M}(2,3)$:
{\footnotesize $$\begin{tabular}{|c|c|c|c|c|}
  \hline
  $f(x)$ & $x^2+2x+2$ & $x^2+x+2$ & $x^2+1$  & $x^2+2x+1$\\
  \hline
  $f_1(x)$ & $x^2+1$ & $x^2+1$   &$(x+1)^2$  &   $x^2+x+1$\\
  \hline
  $\a$    &    $8$    &  $8$ & $4$ & $6$ \\
  \hline
  $\b$    &    $1$    &  $1$ & $1$ & $2$ \\
  \hline
  $g$    &    $10$    &  $10$ & $1$ & $1$ \\
  \hline
\end{tabular}$$}
\vskip 3mm
 (3)\,  $\s\in {\bf M}(3,2)$:
{\footnotesize $$\begin{tabular}{|c|c|c|c|c|}
  \hline
  $f(x)$ & $x^3+1$ & $x^3+x^2+1$ & $x^3+x+1$ & $x^3+x^2+x+1$\\
  \hline
  $\a$    &    $3$    &  $7$ & $7$ & $4$ \\
  \hline
  $\b$    &    $\frac12$    &  $1$ & $1$ & $1$ \\
  \hline
  $g$    &    $1$    &  $7$ & $7$ & $1$ \\
  \hline
\end{tabular}$$}

\vskip 3mm
(4)\, $\s\in {\bf M}(3,3)$:
{\footnotesize $$\begin{tabular}{|c|c|c|c|c|c|}
  \hline
  $f(x)$   & $x^3+1$  & $x^3+2x^2+1$      & $x^3+2x+1$     & $x^3+x^2+2x+1$ & $x^3+2x^2+x+1$\\
  \hline
  $f_1(x)$ & $x^3+2$  & $x^3+2x^2+2x+2$   & $x^3+x^2+x+2$  &   $x^3+2x+2$   &  $x^3+x^2+2$\\
  \hline
  $\a$    &    $6$    &  $26$             & $26$           & $26$           & $26$\\
    \hline
  $\b$    & $\frac23$ &  $2$              & $2$            & $2$            & $2$\\
    \hline
  $g$    &    $19$    &  $136$            & $136$          & $136$          & $136$\\
  \hline
\end{tabular}
$$}
\end{exam}

\end{document}